\documentclass{article}

\usepackage{graphicx}
\usepackage{amssymb,amsfonts,amsmath,amsthm}
\DeclareMathOperator{\spec}{spec}
\DeclareMathOperator{\diag}{diag}
\DeclareMathOperator{\QR}{QR}
\DeclareMathOperator{\Jac}{Jacobi}
\DeclareMathOperator{\CG}{CG}

\DeclareMathOperator{\GMRES}{GMRES}
\DeclareMathOperator{\GOE}{GOE}
\DeclareMathOperator{\GUE}{GUE}
\DeclareMathOperator{\QUE}{QUE}
\DeclareMathOperator{\Genetic}{Genetic}
\DeclareMathOperator{\COSH}{COSH}

\newtheorem{remark}{Remark}

\title{Universality in Numerical Computations with Random Data.  Case Studies.}

\author{Percy Deift\footnote{Courant Institute, New York,
NY, USA},
Govind Menon\footnote{Brown University, Providence, RI, USA},
Sheehan Olver\footnote{The University of Sydney, Sydney, NSW, Australia} and
Thomas Trogdon$^*$}

\begin{document}

\maketitle
\begin{abstract}

The authors present evidence for universality in numerical computations with random data.  Given a (possibly stochastic) numerical algorithm with random input data, the time (or number of iterations) to convergence (within a given tolerance) is a random variable, called the halting time.  Two-component universality is observed for the fluctuations of the halting time, \emph{i.e.}, the histogram for the halting times, centered by the sample average and scaled by the sample variance (see Eqs. \ref{fluctuation} and \ref{halting-fluct} below), collapses to a universal curve, independent of the input data distribution, as the dimension increases.  Thus, up to two components, the sample average and the sample variance, the statistics for the halting time are universally prescribed.  The case studies include six standard numerical algorithms, as well as a model of neural computation and decision making.  { A link to relevant software is provided in \cite{TrogdonNU} for the reader who would like to do computations of his'r own.}

 % The authors perform an empirical study of the distribution of the halting time for many choices of well-known numerical algorithms and well-studied random inputs.  Two-component universality is observed for the halting time:  if a sample is normalized to mean zero and variance one then the empirical histogram for the halting time evidently collapses to a universal curve, independent of the input data distribution, as the dimension increases.  The first six examples that are presented have clear implications for numerical analysis.  The final example of this phenomenon is in a model of neural computation and decision making.

\end{abstract}
\setcounter{section}{1}
%\keywords{universality | random matrix theory | iterative algorithms | computation}

%\abbreviations{iid, independent and identically distributed; CG, conjugate gradient; GMRES, generalized minimal residual; GOE, Gaussian Orthogonal Ensemble; GUE, Gaussian Unitary Ensemble; BE, Bernoulli Ensemble; CS, critically scaled; cLOE, CS Laguerre Orthogonal Ensemble; cPBE, CS positive-definite Bernoulli Ensemble; cSBE, CS shifted Bernoulli Ensemble; cSGE, CS shifted Ginibre Ensemble; QUE, Quartic Unitary Ensemble; BDE, Bernoulli Dirichlet Ensemble; UDE, Uniform Dirichlet Ensemble; COSH, Cosh Unitary Ensemble}

In earlier work \cite{DiagonalRMT}, two of the authors (P.D.~and G.M., together with C.~Pfrang) considered the problem of computing the eigenvalues of a real, $n\times n$ random symmetric matrix $M = (M_{ij})$.  They considered matrices chosen from different ensembles $E$ using a variety of different algorithms $A$.  Let $S_n$ denote the space of real, $n\times n$ symmetric matrices.  Standard eigenvalue algorithms involve iterations of isospectral maps $\varphi = \varphi_A: S_n \rightarrow S_n$, $\spec(\varphi_A(M)) = \spec(M)$ for $M \in S_n$.  If $M \in S_n$ is given, one considers the sequence of matrices $M_{k+1} = \varphi(M_k)$, $k \geq 0$, with $M_0 = M$.  Clearly, $\spec(M_{k+1}) = \spec(M_{k}) = \cdots = \spec(M)$, and under appropriate conditions $M_k = \varphi_A^{(k)}(M)$ converges to a diagonal matrix, $\Lambda = \diag(\lambda_1, \ldots, \lambda _n)$.  Necessarily, the $\lambda_i$'s are the desired eigenvalues of $M$.

In \cite{DiagonalRMT}, the authors discovered the following phenomenon:  For a given accuracy $\epsilon$, a given matrix size $n$ ($\epsilon$ small, $n$ large, in an appropriate scaling range) and a given algorithm $A$, the \emph{fluctuations} in the time to compute the eigenvalues to accuracy $\epsilon$ with the given algorithm $A$, were \emph{universal}, independent of the choice of ensemble $E$.  More precisely, they considered fluctuations in the \emph{deflation time} $T$ (The notion of deflation time is generalized to the notion of \emph{halting time} in subsequent calculations).  Recall that if an $n\times n$ matrix has block form
\begin{align*}
M = \left( \begin{array}{cc} M_{11} & M_{12} \\ M_{21} & M_{22} \end{array} \right)
\end{align*}
where $M_{11}$ is $k \times k$ and $M_{22}$ is $(n-k) \times (n-k)$ for some $1 \leq k \leq n-1$ then one says that the block diagonal matrix $\hat M = \diag( M_{11}, M_{22})$ is \emph{obtained from $M$ by deflation}.  If $\|M_{12}\| = \|M_{21}\| \leq \epsilon$, then the eigenvalues $\{\lambda_i\}$ of $M$ differ from the eigenvalues $\{\hat \lambda_i\}$ of $\hat M$ by $\mathcal O(\epsilon)$.  Let $T = T_{\epsilon,n,A,E}(M)$ be the time ($=$ \# of steps = \# iterations of $\varphi_A$) it takes to deflate a random matrix $M$, chosen from an ensemble $E$, to order $\epsilon$, using algorithm $A$, \emph{i.e.} $T$ is the smallest time such that  for some $k$, $1 \leq k \leq n-1$, $\|(\varphi_A^{(T)}(M))_{12}\| =\|(\varphi_A^{(T)}(M))_{21}\| \leq \epsilon$.   

As explained in \cite{DiagonalRMT}, $T$ is a useful measure of the time required to compute the eigenvalues of $M$:  {Generically,} at worst $\mathcal O(n)$ deflations are needed to compute the eigenvalues of $M$, and at best, $\mathcal O(\log n)$.  The fluctuations $\tau_{\epsilon,n,A,E}(M)$ of $T$ are defined by
\begin{align}\label{fluctuation}
\tau_{\epsilon,n,A,E}(M) = \frac{ T_{\epsilon,n,A,E}(M) - \langle T_{\epsilon,n,A,E} \rangle}{\sigma_{\epsilon,n,A,E}},
\end{align}
where $\langle T_{\epsilon,n,A,E}\rangle$ is the sample average of $T_{\epsilon,n,A,E}(M)$ taken over matrices $M$ from $E$, and $\sigma^2_{\epsilon,n,A,E}$ is the sample variance.  For a given $E$, a typical sample size in \cite{DiagonalRMT} was of order $5,\!000$ to $10,\!000$ matrices $M$, and the output of the calculations in \cite{DiagonalRMT} was recorded in the form of a histogram for $\tau_{\epsilon,n,A,E}$.

\begin{figure}
\centerline{\includegraphics[width=\textwidth]{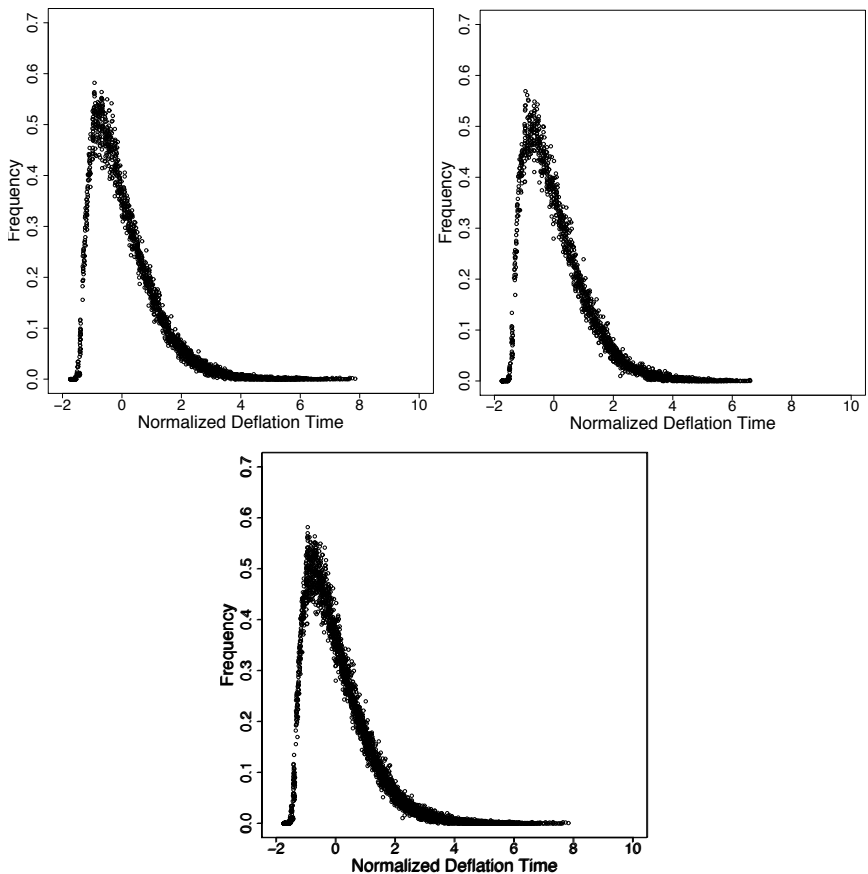}}
\caption{\label{TodaPfrang}  The observation of two-component universality for $\tau_{\epsilon,n,A,E}$ when $A = \QR$. This figure is taken from \cite{DiagonalRMT}.  Overlayed histograms demonstrate the collapse of the histogram of $\tau_{\epsilon,n,A,E}$ to a single curve.  See the Appendix for the definitions of our choices for $E$. In the top-left figure, $E = \GOE$, and 40 histograms for $\tau_{\epsilon,\epsilon,A,E}$, are plotted one on top of the other for $\epsilon = 10^{-k}$, $ k = 2,\!4,\!6,\!8$ and $n = 10,\! 30,\ldots,\! 190$.  The histograms are created with $\approx 10,\!000$ samples.  The top-right figure displays the same information as that in the top-left position, but now for $E = BE$.  In the lower figure, all $40+40$ histograms are overlayed and universality is evident:  the data appears to follow a universal law for the fluctuations. }
\end{figure}

Most of the calculations in \cite{DiagonalRMT} concerned three eigenvalue algorithms: the QR algorithm, the QR algorithm with shifts (the version of QR used in practice), and the Toda algorithm.  The \emph{QR algorithm} is based on the factorization of a(n invertible) matrix $M$ as $M = QR$, where $Q$ is orthogonal and $R$ is upper-triangular with $R_{ii} > 0$.  Given $M \in S_n$, with $M = QR$,  $M' = \varphi_A(M) = \varphi_{\QR}(M) \equiv RQ$.  Clearly, $M' = Q^T M Q \in S_n$ and $\spec (M') = \spec(M)$.
{Practical implementation of the QR algorithm requires the use of a shift, \emph{i.e.} the \emph{QR algorithm with shifts} \cite{Parlett1998}. As shown in \cite{DiagonalRMT}, shifting does not affect universality.}  The \emph{Toda algorithm} involves the solution $M(t)$ of the \emph{Toda equation} $\frac{dM}{dt} = [B(M),M] = B(M) M - M B(M)$, where $B(M) = M_+ - M_+^T$, $M_+$ is the upper triangular part of $M$, and $M(t = 0) = M$.  For all $t > 0$, $\spec(M(t)) = \spec(M)$, and as $t \rightarrow \infty$, we again have $M(t) \rightarrow \Lambda = \diag(\lambda_1, \ldots, \lambda_n)$ where $\{\lambda_i\}$ are the eigenvalues of $M$.  For the convenience of the reader, in Figure~\ref{TodaPfrang}, we reproduce, in particular, histograms for $\tau_{\epsilon,n,A,E}$, from \cite{DiagonalRMT} for the QR algorithm ($A = \QR$) with two different ensembles and varying values of $n$ and $\epsilon$.

From Figure~\ref{TodaPfrang}, we see that eigenvalue computation with the QR algorithm exhibits two-component universality, \emph{i.e.}, the fluctuations $\tau_{\epsilon,n,A,E}$ obey a universal law for all ensembles $E$ under consideration.  The same is true for all three algorithms considered in \cite{DiagonalRMT}:  The laws are different, however, for different algorithms $A$.

In the current paper, the work in \cite{DiagonalRMT} has been extended in various ways as follows.   All matrix ensembles are described in the Appendix.

\subsection{The Jacobi Algorithm}\label{sec:Jacobi}

\begin{figure}
\centerline{\includegraphics[width=\textwidth]{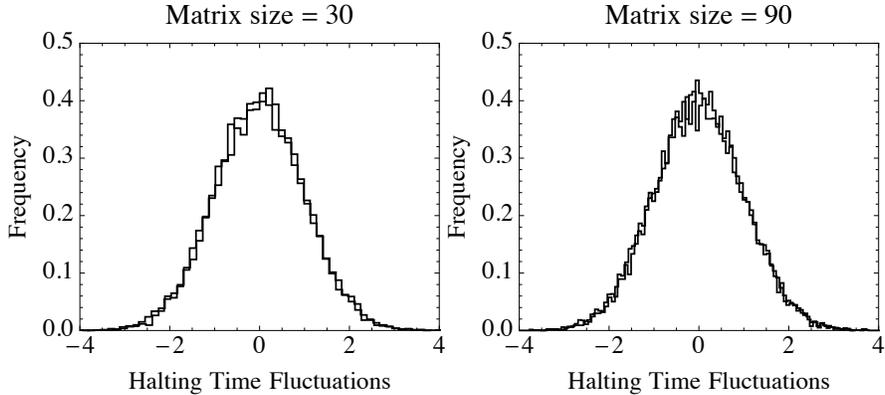}}
\caption{\label{Jacobi}  The observation of two-component universality for $\tau_{\epsilon,n,A,E}$ when $A = \Jac$, $E = \GOE,~\mathrm{BE}$ and $\epsilon = \sqrt{n}\, 10^{-10}$. The left figure displays two histograms, one on top of the other, one for $\GOE$ and one for $\mathrm{BE}$, when $n = 30$.  The right figure displays the same information for $n = 90$.  All histograms are produced with $16,\!000$ samples. We see two-component universality emerge for $n$ sufficiently large: the histograms follow a universal (independent of $E$) law.}
\end{figure}

In the first set of computations, the authors consider the eigenvalue problem for random matrices $M \in S_n$ using the Jacobi algorithm (see, \emph{e.g.} \cite{Golub2013}): for $M \in S_n$, choose $i < j$ such that $|M_{ij}| \geq \max_{1 \leq i' < j' \leq n} |M_{i'j'}|$, and let $G^{(ij)} \equiv G^{(ij)}(\theta)$ be the corresponding \emph{Givens rotation matrix}:  $G^{(ij)}_{i'j'} = \delta_{i'j'}$, for $i',j' \neq i,j$, and
\begin{align*}
\left[ \begin{array}{cc} G^{(ij)}_{ii} & G^{(ij)}_{ij} \\ G^{(ij)}_{ji} & G^{(ij)}_{jj} \end{array} \right] = \left[ \begin{array}{cc} \cos(\theta) & \sin(\theta) \\ - \sin (\theta) & \cos(\theta) \end{array} \right], ~~ (G^{(ij)})^T G^{(ij)} = I.
\end{align*}
Here $\theta = \theta(M)$ is chosen so that $((G^{(ij)})^T M G^{(ij)})_{ij} = 0$  and then $\varphi_{\Jac}(M) \equiv (G^{(ij)})^T M G^{(ij)}$. Clearly, $M' = \varphi_{\Jac}(M) \in S_n$ and $\spec(M') = \spec(M)$ and again (see \cite{Golub2013}), $M_k = \varphi^{(k)}_{\Jac}(M) \rightarrow \Lambda = \diag(\lambda_1,\ldots,\lambda_n)$.  The Jacobi algorithm has a very different character from QR-Toda type algorithms which are intimately connected  to completely integrable Hamiltonian systems (see \cite{DeiftEigenvalue} and the references therein)\footnote{The Jacobi algorithm is well-suited to parallel computation, and also has other advantages over QR in the context of modern, large-scale computation ( see \emph{e.g.} \cite{Demmel1992}).}  Deflation, which is a useful measure for eigenvalue computation times for QR/Toda type algorithms, is not useful for the Jacobi algorithm.  In place of $T_{\epsilon,n, A, E}$, we record the \emph{halting time} $k_{\epsilon,n,A,E}$: the number of iterations it takes for the Jacobi algorithm to reduce the Frobenius norm of the off-diagonal elements to be less than a given\footnote{This is sufficient to conclude that one element on the diagonal of the transformed matrix is within $\epsilon n^{-1/2}$ of an exact eigenvalue of the original matrix.} $\epsilon$. Histograms are produced for an appropriate analog of $\tau_{\epsilon,n,A,E}$:
\begin{align}\label{halting-fluct}
\tau_{\epsilon,n,A,E}(M) = \frac{k_{\epsilon,n,A,E}(M) - \langle k_{\epsilon,n,A,E} \rangle}{\sigma_{\epsilon,n,A,E}}.
\end{align}
Computations for $A = \Jac$ are given in Figure~\ref{Jacobi}.  Again, two-component universality is evident.

\subsection{Ensembles with Dependent Entries}\label{sec:Depend}

\begin{figure}
\centerline{\includegraphics[width=\textwidth]{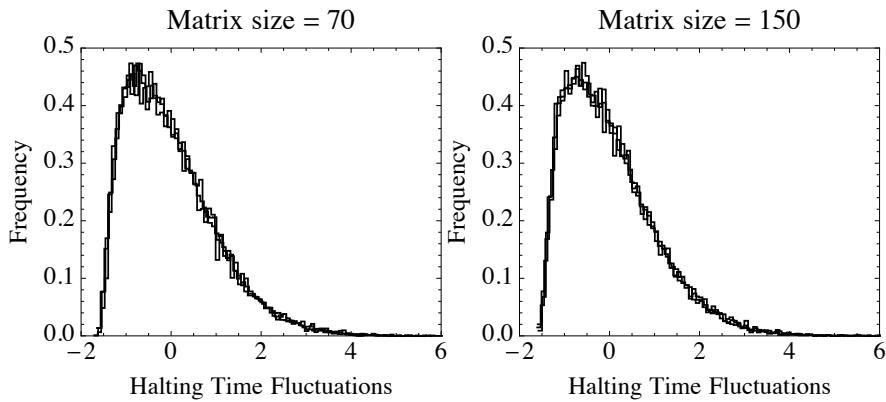}}
\caption{\label{QR}  The observation of two-component universality for $\tau_{\epsilon,n,A,E}$ when $A = \QR$, $E = \mathrm{QUE}, ~\COSH, ~\mathrm{GUE}$ and $\epsilon = 10^{-10}$.  Here we are using deflation time ( = halting time), as in \cite{DiagonalRMT}. The left figure displays three histograms, one each for $\GUE,~\COSH$ and $\QUE$, when $n = 70$.  The right figure displays the same information for $n = 150$.  All histograms are produced with $16,\!000$ samples.  Two-component universality emerges for $n$ sufficiently large: the histograms follow a universal (independent of $E$) law.  This is surprising because $\COSH$ and $\QUE$ have eigenvalue distributions that differ significantly from $\GUE$ in that they do not follow the so-called \emph{semi-circle law}. {These histograms appear to collapse to the same curve in Figure~\ref{TodaPfrang}.}  This is a further surprise, given the well-known fact that Orthogonal and Unitary Ensembles give rise to different (eigenvalue) universality classes.}
\end{figure}

In all the above cases, including the calculations for the Jacobi algorithm, the matrices $M$ are real and the entries $M_{ij}$ are independent, subject only to the symmetry requirement $M_{ij} = M_{ji}$.  In the second set of computations in the present paper, the authors consider $n \times n$ Hermitian matrices $M = M^*$ taken from various unitary  ensembles (see \emph{e.g} \cite{MehtaRM}) with probability distributions proportional to $e^{-n \mathrm{tr} V(M)} dM$ where $V: \mathbb R \rightarrow \mathbb R$ grows sufficiently rapidly as $|x| \rightarrow \infty$, and $dM$ is Lebesgue measure on the algebraically independent entries $M_{ij} = \mathrm{Re}\, M_{ij} + \sqrt{-1}\, \mathrm{Im}\, M_{ij}$ of $M$.  Unless $V(x)$ is proportional to $x^2$, the entries of $M$ for such ensembles are dependent, and it is a non-trivial matter to sample the matrices.  A novel technique for sampling such unitary  ensembles was introduced recently \cite{Olver2014} by two of the authors, S.O. and T.T., together with N. R. Rao, taking advantage of the representation of the eigenvalues of $M$ as a determinantal point process whose kernel is given in terms of orthogonal polynomials (see also \cite{Li2013}).  Using this sampling technique, the authors of the present paper have considered the QR algorithm for various unitary  ensembles\footnote{Here $M = QR$ where $Q$ is unitary and again $R$ is upper triangular with $R_{ii} > 0$.}.  Histograms for the halting (= deflation) time fluctuations $\tau_{\epsilon,n,A,E}$, $A = QR$, are given in Figure~\ref{QR} and again two-component universality is evident.

\subsection{The Conjugate Gradient Algorithm}\label{sec:CG}

\begin{figure}
\centerline{\includegraphics[width=\textwidth]{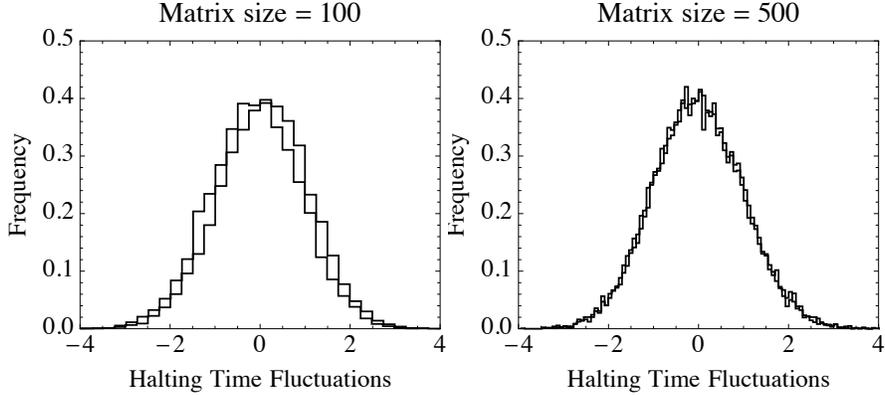}}
\caption{\label{CG} The observation of two-component universality for $\tau_{\epsilon,n,A,E}$ when $A = \CG$ and $E = \mathrm{cLOE}, ~\mathrm{cPBE}$ with $\epsilon = 10^{-10}$.  The left figure displays two histograms, one for $\mathrm{cLOE}$ and $\mathrm{cPBE}$, when $n = 100$.  The right figure displays the same information for $n = 500$.  All histograms are produced with $16,\!000$ samples.  Two-component universality is evident for $n$ sufficiently large: the histograms follow a universal (independent of $E$) law.   The critical scaling (see Appendix~\ref{app:pd}) has significant impact on the distribution of the condition number and forces $\langle \tau_{\epsilon,n,A,E} \rangle \approx n \alpha$, $\alpha < 1$.  If the scaling $m = 2n$  is chosen in the ensemble $E$ then the CG method converges too quickly and the halting time tends to take only 10-15 different values for each value of $m$. No interesting limiting statistics are present.  Conversely, if $m = n$ the CG method converges slowly ($\langle k_{\epsilon,m,A,E} \rangle \gg m$) and rounding errors dominate the computation.  Experiments do not indicate two-component universality if $m = 2n$ or $m = n$.  The scaling $m = n + 2 \lfloor \sqrt{n} \rfloor$ identifies a critical scaling region. Within this scaling region,  we see two-component universality emerge for $n$ sufficiently large: the histograms follow a universal (independent of $E$) law.  }
\end{figure}

In a third set of computations in this paper, the authors start to address the question of whether two-component universality is just a feature of eigenvalue computation, or is present more generally in numerical computation.  In particular, the authors consider the solution of the linear system of equations $Wx=b$ where $W$ is real and positive definite, using the conjugate gradient (CG) method.  The method is iterative (see \emph{e.g.} \cite{Saad2003} and also Remark~\ref{rmk:scaling} below) and at iteration $k$ of the algorithm an approximate solution $x_k$ of $Wx=b$ is found and the residual $r_k = Wx_k-b$ is computed.  For any given $\epsilon > 0$, the method is halted when\footnote{The notation $\|\cdot\|_2$ is used to denote the standard $\ell^2$ norm on $n$-dimensional Euclidean space} $\|r_k\|_2 < \epsilon$, and the halting time $k_{\epsilon}(W,b)$ recorded.  The authors consider $n \times n$ matrices $W$ chosen from two different positive definite ensembles $E$ (see Appendix~\ref{app:pd}) and vectors $b = (b_j)$ chosen independently with iid~entries $\{b_j\}$.  Given $\epsilon$ (small) and $n$ (large), and $(W,b) \in E$, the authors record the halting time $k_{\epsilon,n,A,E}$, $A = \CG$, and compute the fluctuations $\tau_{\epsilon,n,A,E}(W,b)$. The histograms for $\tau_{\epsilon,n,A,E}$ are given in Figure~\ref{CG}, and again, two-component universality is evident.

\subsection{The GMRES Algorithm}\label{sec:GMRES-un}

\begin{figure}
\centerline{\includegraphics[width=\textwidth]{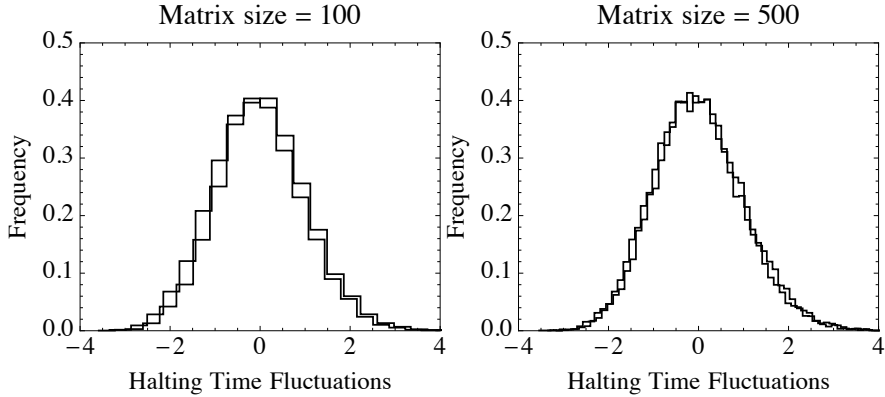}}
\caption{\label{GMRES-un} The observation of two-component universality for $\tau_{\epsilon,n,A,E}$ when $A = \GMRES$, $E = \mathrm{cSGE}, ~\mathrm{cSBE}$ and $\epsilon = 10^{-8}$. The left figure displays two histograms, one for $\mathrm{cSGE}$ and one for $\mathrm{cSBE}$, when $n = 100$.  The right figure displays the same information for $n = 500$.  All histograms are produced with $16,\!000$ samples.   The critically scaled ensembles cSBE and cSGE are of the form $I + X_n$ with $ \|X_n\| \approx 2$. If the matrix is too close to the identity, the halting time will take almost constant values, \emph{i.e.} $k_{\epsilon,n,A,E} = 8$, independent of $n$.  If the matrix is too far from the identity, the fact that it is unstructured makes GMRES perform poorly and the algorithm typically completes in $n$ steps, the maximum possible number of iterations (see Remark~\ref{rmk:scaling} below).  With the proper scaling of $X$, we see two-component universality emerge for $n$ sufficiently large: the histograms follow a universal (independent of $E$) law.}
\end{figure}

In a fourth set of computations, the authors again consider the solution of $Wx =b$ but here $W$ has the form $I + X$ and $X \equiv X_n$ is a random, real non-symmetric matrix and $b= (b_j)$ is independent with uniform iid~entries $\{b_j\}$.  As $W = I +X$ is (almost surely) no longer  positive definite the conjugate gradient algorithm breaks down, and the authors solve $(I+X)x = b$ using the Generalized Minimal Residual (GMRES) algorithm \cite{GMRES-original}.  Again, the algorithm is iterative and at iteration $k$ of the algorithm an approximate solution $x_k$ of $(I+X)x = b$ is found and the residual $r_k = (I+X)x_k -b$ is computed.  As before, for any given $\epsilon > 0$, the method is halted when $\|r_k\|_2 < \epsilon$ and $k_{\epsilon,n,A,E}(X,b)$ is recorded.  As in the conjugate gradient problem (Section~\ref{sec:CG}), the authors compute the histograms for the fluctuations of the halting time $\tau_{\epsilon,n,A,E}$ \eqref{halting-fluct} for two ensembles $E$, where now $A = \GMRES$.  The results are given in Figure~\ref{GMRES-un}, where again two-component universality is evident.

\begin{remark}\label{rmk:scaling}
The computations in Sections~\ref{sec:CG} and \ref{sec:GMRES-un} are particularly revealing for the following reason.  Both the CG and GMRES algorithms proceed by generating approximations $x_n$ to the solution in progressively larger subspaces $V_k$ of $\mathbb R^n$, $x_k \in V_k$, $\dim V_k = k$ (almost surely).  These algorithms terminate in at most $n$ steps, in the absence of rounding errors.  If the matrix $W$ in the case of CG, or $I +X$ in the case of GMRES, is too close to the identity, then the algorithm will converge in $\mathcal O(1)$ steps, essentially independent of $n$.   On the other hand, if $W$ or $I + X$ is too far from the identity, the algorithm will converge only after $n$ steps (GMRES) or be dominated by rounding errors (CG).  Thus in both cases there are no meaningful statistics.  What the calculations in Sections~\ref{sec:CG} and \ref{sec:GMRES-un} reveal is that if the ensembles for CG and GMRES are such that the matrices $W$ and $I + X$, respectively, are typically not too close to, and not too far, from the identity, then the algorithms exhibit significant statistical fluctuations, and two-component universality is immediately evident. (for more discussion see the captions for Figures~\ref{CG} and \ref{GMRES-un}).  Analogous considerations apply in Section~\ref{sec:GMRES-Dir} below.
\end{remark}

\subsection{Discretization of a Random PDE}\label{sec:GMRES-Dir}

\begin{figure}
\centerline{\includegraphics[width=\textwidth]{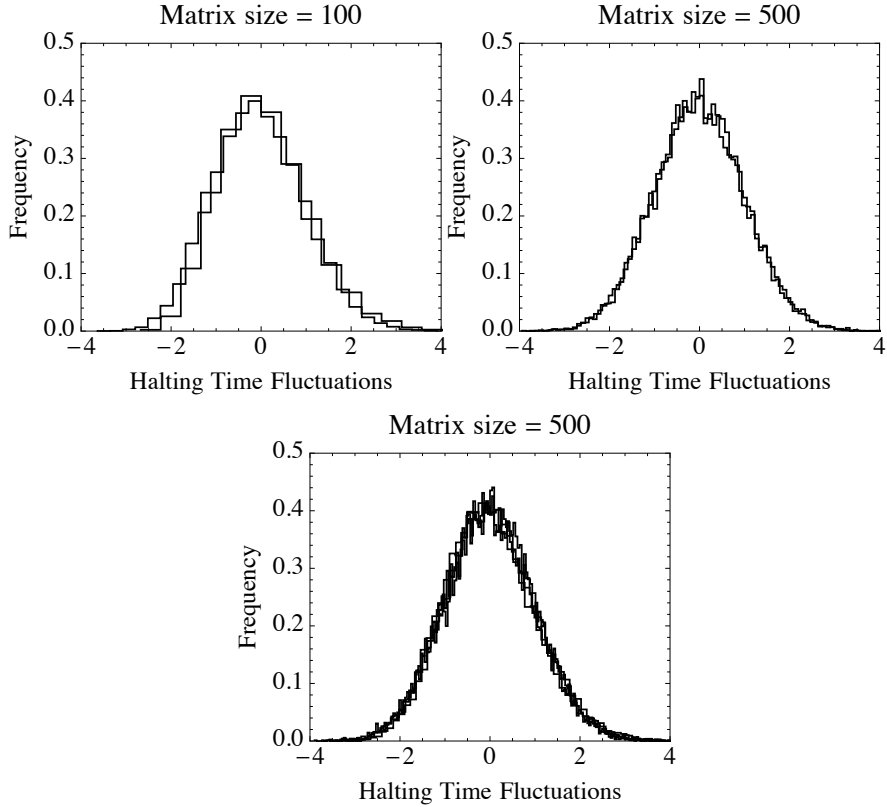}}
\caption{\label{GMRES-Dir}  The observation of two-component universality for $\tau_{\epsilon,n,A,E}$ when $A = \GMRES$, $E = \mathrm{UDE}, ~\mathrm{BDE}$ and $\epsilon = 10^{-8}$. The left figure displays two histograms, one for $\mathrm{UDE}$ and one for $\mathrm{BDE}$, when $n = 100$.  The right figure displays the same information for $n = 500$.  The bottom figure consists of four histograms, two taken from Figure~\ref{GMRES-un} ($E = \mathrm{cSGE},~\mathrm{cSBE}$) and two taken from the right figure above ($E = \mathrm{UDE},~\mathrm{BDE}$).  All histograms are produced with $16,\!000$ samples.  It is interesting to note two properties.   First, as we observe from our computations, BDE and UDE are of the form $I + X_n$ where $X_n$ has a norm that grows proportional to some fractional power of $n$.  While this type of growth in the case of Section~\ref{sec:GMRES-un} (Figure~\ref{GMRES-un}) would cause GMRES to take its maximum possible number of iterations, that is $k = n$, nevertheless, in the context of Section~\ref{sec:GMRES-Dir}, non-trivial statistics emerge.  In light of Remark~\ref{rmk:scaling}, we conjecture that structure is necessary for GMRES to perform well when the perturbation of the identity has an unbounded spectral radius in the large $n$ limit.  The second, and most important feature, is that two-component universality for matrices of the form $I + X_n$ persists as the computations are moved from structured randomness (UDE and BDE) to unstructured randomness (cSBE and cSGE): the histograms follow a universal (independent of $E$) law.}
\end{figure}

In a fifth set of computations, the authors raise the issue of whether two-component universality is just a feature of finite-dimensional computation, or is also present in problems which are intrinsically infinite dimensional.  In particular, is the universality present in numerical computations for PDEs?  As a case study, the authors consider the numerical solution of the Dirichlet problem $\Delta u = 0$ in  a star-shaped region $\Omega \subset \mathbb R^2$ with $ u = f$ on $\partial \Omega$. The boundary is described by a periodic function of the angle $\theta$, $r = r(\theta)$, and similarly $f = f(\theta)$, $0 \leq \theta \leq 2 \pi$.  Two ensembles, BDE and UDE (as described in Appendix~\ref{app:Dir}), are derived from a discretization of the problem with specific choices for $r$, defined by a random Fourier series.  The boundary condition $f$ is chosen randomly by letting $\{f(\frac{2 \pi j}{n})\}_{j=0}^{n-1}$ be iid~uniform on $[-1,1]$.  Histograms for the halting time $\tau_{\epsilon,n,A,E}$ from these computations are given in Figure~\ref{GMRES-Dir} and again, two-component universality is evident.  What is surprising, and quite remarkable, about these computations is that the histograms for $\tau_{\epsilon,500,A,E}$ in this case are the \emph{same} as the histograms for $\tau_{\epsilon,500,A,E}$ in Figure~\ref{GMRES-un} (see Figure~\ref{GMRES-Dir} for the overlayed histograms).  In other words, UDE and BDE are structured with random components, whereas cSGE and cSBE have no structure, yet they produce the same statistics (modulo two components).

\subsection{A Genetic Algorithm}\label{sec:Genetic}

In all the computations discussed so far, the randomness in the computations\footnote{Aside from round-off errors, see comments below Figure~\ref{CG}} resides in the initial data.  In the sixth set of computations, the authors consider an algorithm which is intrinsically stochastic.  They consider a genetic algorithm to compute \emph{Fekete points} (see \cite[p.~142]{SaffPotential}).  Such points $P^* = (P_1^*,P_2^*, \ldots, P^*_N) \in \mathbb R^N$ are the global minimizers of the objective function
\begin{align*}
H(P) = \frac{2}{N(N-1)} \sum_{1 \leq i \neq j \leq N} \log |P_i- P_j|^{-1} + \frac{1}{N} \sum_{i=1}^N V(P_i)
\end{align*}
for real-valued functions $V = V(x)$ which grow sufficiently rapidly as $|x| \rightarrow \infty$.  It is well-known (see, \emph{e.g.} \cite{SaffPotential})  that  as $N \rightarrow \infty$, the counting measures $\delta_{P^*} = \frac{1}{N} \sum_{i=1}^N \delta_{P_i^*}$ converge to the so-called equilibrium measure $\mu_V$ which plays a key role in the asymptotic theory of the orthogonal polynomials generated by measure $e^{-NV(x)}dx$ on $\mathbb R$.  Genetic algorithms involve two basic components , ``mutation'' and ``crossover''.  The authors implement the genetic algorithm in the following way.

\paragraph{The Algorithm}

Fix a distribution $\mathfrak D$ on $\mathbb R$.  Draw an initial population $\mathcal P_0 = \mathcal P = \{P_i\}_{i=1}^n$ consisting of $n = 100$ vectors in $\mathbb R^N$, $N$ large, with elements that are iid uniform on $[-4,4]$.  The random map $F_{\mathfrak D}(\mathcal P): (\mathbb R^N)^n \rightarrow (\mathbb R^N)^n$ is defined by one of the following two procedures:
\begin{itemize}
\item \textbf{Mutation}: Pick one individual $P \in \mathcal P$ at random (uniformly).  Then pick two integers $n_1,~n_2$ from $\{1,2,\ldots,N\}$ at random (uniformly and independent).  Three new individuals are created.
\begin{itemize}
\item $\tilde P_1$ --- draw $n_1$ iid numbers $\{x_1, \ldots, x_{n_1} \}$ from $\mathfrak D$ and perturb the first $n_1$ elements of $P$ : $(\tilde P_1)_i = (P)_i + x_i$, $i = 1, \ldots, n_1$, and $(\tilde P_1)_i = (P)_i$ for $i > n_1$.
\item $\tilde P_2$ --- draw $N - n_2$ iid numbers $\{y_{n_2+1},\ldots,y_{N}\}$ from $\mathfrak D$ and perturb the last $N- n_2$ elements of $P$: $(\tilde P_2)_i = (P)_i + y_i$, $i = n_2+1, \ldots, N$, and $(\tilde P_2)_i = (P)_i$ for $i \leq n_2$.
\item $\tilde P_3$ --- draw $|n_1-n_2|$ iid numbers $\{z_{1},\ldots,z_{|n_1-n_2|}\}$ from $\mathfrak D$ and perturb elements $n_1^*=1+\min(n_1,n_2)$ through $n_2^*=\max(n_1,n_2)$: $(\tilde P_3)_{i} = (P)_i + z_{i-n_1^*+1}$, $i = n_1^*, \ldots, n_2^*$, and $(\tilde P_3)_i = (P)_i$ for $i \not\in \{n_1^*,\ldots,n_2^*\}$.

\end{itemize}

\item \textbf{Crossover}:  Pick two individuals $P,~Q$ from $\mathcal P$ at random (independent and uniformly).  Then pick two numbers $n_1,~n_2$ from $\{1,2,\ldots,N\}$ (independent and uniformly).  Two new individuals are created.

\begin{itemize}
\item $\tilde P_4$ --- Replace the $n_1$th element of $P$ with the $n_2$th element of $Q$ and perturb it (additively) with a sample of $\mathfrak D$.
\item $\tilde P_5$ --- Replace the $n_1$th element of $Q$ with the $n_2$th element of $P$ and perturb it (additively) with a sample of $\mathfrak D$.
\end{itemize}
\end{itemize}
At each step, the application of either crossover or mutation is chosen  with equal probability.  The new individuals are appended\footnote{After mutation we have $\tilde {\mathcal P} = \mathcal P \cup \{\tilde P_1, \tilde P_2, \tilde P_3\}$ and after crossover, $\tilde {\mathcal P} = \mathcal P \cup \{\tilde P_4, \tilde P_5\}$} to $\mathcal P$ and $\mathcal P \mapsto \mathcal P' = F_{\mathfrak D}(\mathcal P) \in (\mathbb R^N)^n$ is constructed by choosing the 100 $P_i$'s in $\tilde {\mathcal P}$ which yield the smallest values of $H(P)$.  The algorithm produces a sequence of populations $\mathcal P_1, \mathcal P_2,\ldots,\mathcal P_k, \ldots$ in $(\mathbb R^N)^n$, $\mathcal P_{k+1} = F_{\mathfrak D}(\mathcal P)$, $n = 100$, and halts, with halting time recorded, for a given $\epsilon$, when $\min_{P \in \mathcal P_k} H(P) - \inf_{P \in \mathbb R^N}  H(P) < \epsilon$.

\begin{figure}
\centerline{\includegraphics[width=\textwidth]{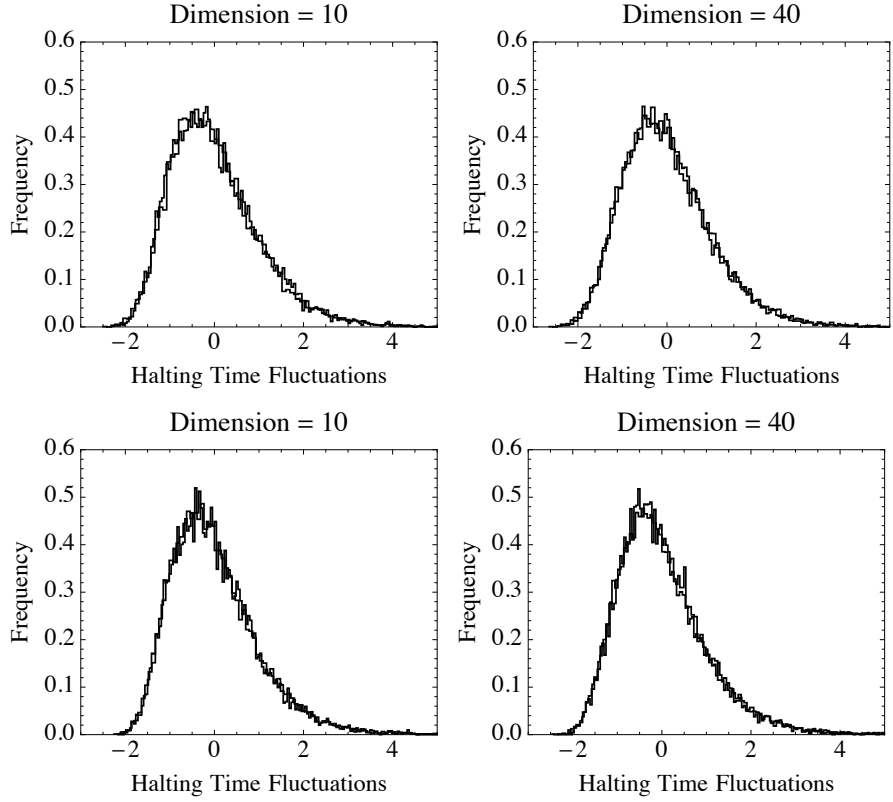}}
\caption{\label{Genetic}  The observation of two-component universality for $\tau_{\epsilon,N,A,E}$ when $A = \Genetic$, $\epsilon = 10^{-2}$ and $E \simeq \mathfrak D$ where $\mathfrak D$ is chosen to be either uniform on $[-1/(10 N), 1/(10 N)]$ or taking values $\pm 1/(10 N)$ with equal probability. The top row is created with the choice $V(x) = x^2$ and the bottom row with $V(x) = x^4-3x^2$.  Each of the plots in the left column displays two histograms, one for each choice of $\mathfrak D$ when $N = 10$.  The right column displays the same information for $N = 40$.   All histograms are produced with $16,\!000$ samples. It is evident that the histograms collapse onto a universal curve, one for each $V$.  }
\end{figure}

The histograms for the fluctuations $\tau_{\epsilon,N,A,E}$, with $A = \Genetic$ are given in Figure~\ref{Genetic}, for two choices of $V$, $V(x) = x^2$ and $V(x) = x^4 - 3 x^2$, and different choices of $E \simeq \mathfrak D$.  For $V(x) = x^2$ $\inf_{P \in \mathbb R^N}  H(P)$ is known explicitly, and for $V(x) = x^4- 3x^2$, $\inf_{P \in \mathbb R^N}  H(P)$ is approximated by a long run of the genetic algorithm.  As before, two-component universality is evident.

\subsection{Curie--Weiss Model}\label{sec:CurieWeiss}

In the seventh and final set of computations, the authors pick up on a common notion in neuroscience that the human brain is a computer with software and hardware.  If this is indeed so, then one may speculate that two-component universality should certainly be present in some cognitive actions.  Indeed, such a phenomenon is in evidence in the recent experiments of Bakhtin and Correll \cite{Bakhtin2012}.  In \cite{Bakhtin2012},  data  from experiments with 45 human participants was analyzed.  The participants are shown 200 pairs of images.  The images in each pair consist of nine black disks of variable size.  The disks in the images within each pair have approximately the same area so that there is no \emph{a priori} bias.  The participants are then asked to decide which of the two images covers a larger (black) area and the time $T$ required to make a decision is recorded.  For each participant, the decision times for the 200 pairs are collected and the fluctuation histogram\footnote{In \cite{Bakhtin2012} the authors do not display the histogram for the fluctuations directly, but such information is easily inferred from their figures (see Figure~6 in \cite{Bakhtin2012}).} is tabulated.  The experimental results are in good agreement with a dynamical Curie-Weiss model frequently used in describing decision processes \cite{Bakhtin2011}.  As each of the 45 participants operates, presumably, in his'r own stochastic neural environment, this is a remarkable demonstration of two-component universality in cognitive action.

At its essence the Curie--Weiss model is Glauber dynamics on the hypercube $\{-1,1\}^N$ with a microscopic approximation of a drift-diffusion process.  Consider $N$ variables $\{X_i(t)\}_{i=1}^N$, $X_i(t) \in \{-1,1\}$. The state of the system at time $t$ is $X(t) = (X_1(t), X_2(t), \ldots, X_N(t))$.  The transition probabilities are given through the expressions
\begin{align*}
\mathbb P(X_i(t+\Delta t) \neq X_i(t) | X(t) = x ) = c_i(x) \Delta t + o(\Delta t),
\end{align*}
where $c_i(x)$ is the spin flip intensity. The observable considered is  $M(X(t)) = \frac{1}{N} \sum_{i=1}^N X_i(t) \in [-1,1]$,  and the initial state of the system is chosen so that $M(X(0)) = 0$, a state with no \emph{a priori} bias, as in the case of the experimental setup. The halting (or decision) time for this model is $k = \inf\{t: |M(X(t))| \geq \epsilon\}$, the time at which the system makes a decision.  Here $\epsilon \in (0,1)$ may not be small.

This model is simulated by first sampling an exponential random variable with mean $\lambda(t) = \left(\sum_i c_i(X(t))\right)^{-1}$ to find the time $\Delta t$ at which the system changes state. Sampling the random variable $Y$, $\mathbb P(Y = i) = c_i(X(t)) \lambda(t)$, $i = 1,2,\ldots,N$ produces an integer $j$, determining which spin flipped. Define $X_i(t+s) \equiv X_i(t)$ if $s \in [0,\Delta t)$ for $i = 1,2,\ldots,N$ and $X_i(t + \Delta t) \equiv X_i(t)$, $X_j(t+ \Delta t) \equiv -X_j(t)$ for $i \neq j$.  This procedure is repeated with $t$ replaced by $t+ \Delta t$ to evolve the system.

\begin{figure}
\centerline{\includegraphics[width=\textwidth]{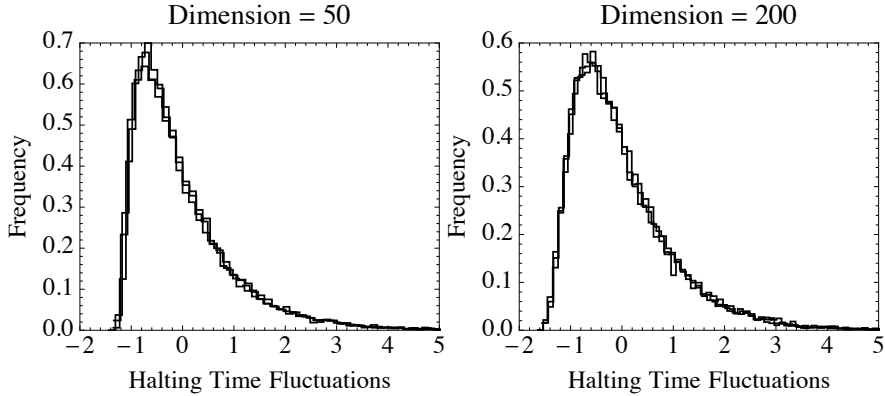}}
\caption{\label{CurieWeiss}
The observation of two-component universality for $\tau_{\epsilon,N,A,E}$ when $A = \text{Curie--Weiss}$, $E \simeq o_i,~ u_i, ~ v_i$, $\epsilon = .5$ and $\beta = 1.3$.  The left figure displays three histograms, one for each choice of $E$ when $N = 50$.  The right figure displays the same information for $N = 200$.   All histograms are produced with $16,\!000$ samples. The histogram for $E = o_i$ corresponds to the case studied in \cite{Bakhtin2012,Bakhtin2011}.  It is clear from these computations that the fluctuations collapse on to the universal curve for $E = o_i$.  Thus, reasonable changes in the spin flip intensity do not appear to change the limiting histogram.  This indicates why the specific choice made in \cite{Bakhtin2012} of $E = o_i$ is perhaps enough to capture the behavior of many individuals. }
\end{figure}

Central to the application of the model is the assumption on the statistics of the spin flip intensity $c_i(x)$.   If one changes the basic statistics of the $c_i$'s, will the limiting histograms for the fluctuations of $k$ be affected as $N$ becomes large?  In response to this question the authors consider the following choices for $E \simeq c_i(x)$ ($\beta = 1.3$): {$c_i(x) = o_i(x) = e^{-\beta x_i M(x)}$ (the case studied in \cite{Bakhtin2012}), $c_i(x) = u_i(x) = e^{-\beta x_i(M(x)-M^3(x)/5)}$, or $c_i(x) = v_i(x) = e^{-\beta x_i(M(x) + M^8(x))}$.}
The resulting histograms for the fluctuations $\tau_{\epsilon,N,A,E}$ of $T$ are given in Figure~\ref{CurieWeiss}. Once again, two-component universality is evident.  Thus the universality in the decision process models mirrors the universality observed among the 45 participants in the experiment of Bakhtin and Correll.

%\section{Discussion}

\section{Conclusions}

{Two distinct themes are combined in this work: (1) the notion of universality in random matrix theory and statistical physics;  (2) the use of random ensembles in scientific computing. The origin of both these ideas dates to the 1950s in the work of (1) Wigner~\cite{MehtaRM,Wigner1951}, and (2) von Neumann and Goldstine~\cite{Goldstine1951}. There has been considerable progress in the rigorous understanding of universality in random matrix theory~(see e.g.~\cite{DeiftRandom4,Erdos2012} and the references therein). In contrast, the performance of numerical algorithms on random ensembles is less understood, though results in this area include probabilistic bounds for condition numbers and halting times for numerical algorithms~\cite{Demmel1988,Edelman1988,Smale1985}.

The work presented here reveals empirical evidence for two-component universality in several numerical algorithms. The results of~\cite{DiagonalRMT} and Sections \ref{sec:Jacobi}-\ref{sec:GMRES-Dir} reveal universal fluctuations of halting times for iterative algorithms in numerical linear algebra on random matrix ensembles with both dependent and independent entries. In each instance, the process of numerical computation on a random matrix may be viewed as the evolution of a random ensemble by a deterministic dynamical system. In a similar light, the algorithms of Section~\ref{sec:Genetic} and \ref{sec:CurieWeiss} may be seen as stochastic dynamical systems with that in Section~\ref{sec:CurieWeiss} having a close connection with neural computation. In all these examples, the empirical observations presented here suggest new  universal phenomena in non-equilibrium statistical mechanics. The results of Section~\ref{sec:GMRES-un} and \ref{sec:GMRES-Dir}
reveal that numerical computations with a structured ensemble with some random components may have the same statistics (modulo two-components) as an unstructured ensemble. This brings to mind the situation in the 1950s when Wigner introduced random matrices as a model for scattering resonances of neutrons off heavy nuclei: the neutron-nucleus system has a well-defined and structured Hamiltonian, but nevertheless the resonances for neutron scattering are well-described statistically by the eigenvalues of an (unstructured) random matrix.}

\section*{Materials}
All algorithms discussed here are implemented in \texttt{Mathematica}.  A package is available for download \cite{TrogdonNU} that contains all relevant data and the code to generate this data.  The package supports parallel evaluation for most algorithms and runs easily on personal computers.

\appendix
\setcounter{section}{1}

\subsection{Gaussian Ensembles}\label{app:gauss}

The Gaussian Orthogonal Ensemble (GOE) is given by  $(X+X^T)/\sqrt{4n}$ where $X$ is an $n \times n$ matrix of {standard} iid~Gaussian variables.  The Gaussian Unitary Ensemble (GUE) is given by $(X+X^*)/\sqrt{8n}$ where $X$ is an $n \times n$ matrix of {standard} iid~complex Gaussian variables.

\subsection{Bernoulli Ensemble}\label{app:bernoulli}

The Bernoulli Ensemble (BE) is given by an $n \times n$ matrix $X$ consisting of iid~random variables that take the values $\pm 1/\sqrt{n}$ with equal probability subject only to the constraint $X^T = X$.

\subsection{Positive Definite Ensembles}\label{app:pd}

The critically-scaled Laguerre Orthogonal Ensemble (cLOE) is given by $ W = XX^T/m$ where $X$ is an $n \times m$ matrix with {standard} iid Gaussian entries. The critically-scaled positive definite Bernoulli ensemble (cPBE) is given by $W = XX^T/m$ where $X$ is an $n \times m$ matrix consisting of iid~Bernoulli variables taking the values $\pm 1$ with equal probability. {In both cases,} the critical scaling refers to the choice $m = n + 2\lfloor \sqrt{n} \rfloor$.

\subsection{Shifted Ensembles}\label{app:shift}

The critically-scaled shifted Bernoulli Ensemble (cSBE) is given by $I + X/\sqrt{n}$ where $X$ is an $n \times n$ matrix consisting of iid~Bernoulli variables taking the values $\pm 1$ with equal probability.  The critically-scaled shifted Ginibre Ensemble (cSGE) is given by $I + X/\sqrt{n}$ where $X$ is an $n \times n$ matrix of {standard} iid~Gaussian variables. With this scaling $ \mathbb P( |\|X/\sqrt{n}\| -2| > \epsilon) \rightarrow 0$ as $n \rightarrow \infty$ \cite{Geman1980}.

\subsection{Unitary Ensembles}\label{app:unitary}

The Quartic Unitary Ensemble (QUE) is a complex, unitary ensemble with probability distribution proportional to $e^{-n \mathrm{tr} M^4}dM$.  The Cosh Unitary Ensemble (COSH) has its distribution proportional to $e^{- \mathrm{tr} \cosh M}dM$.

\subsection{Dirichlet Ensembles}\label{app:Dir}
We consider the numerical solution of the equation $\Delta u = 0$ in $\Omega$ and $u = f$ on $\partial \Omega$.  Here we let $\Omega$ be the star-shaped region interior to the curve $(x,y) = (r(\theta) \cos(\theta),r(\theta) \sin(\theta))$  where $r(\theta)$  for $ 0 \leq \theta < 2 \pi$ is given by  $r(\theta) = 1 + \sum_{j=1}^m (X_j \cos(j \theta) + Y_j \sin(j \theta)),$
and $X_j$ and $Y_j$ are iid~random  variables {on} $[-1/(2m),1/(2m)]$.   The boundary integral equation
\begin{align*}
 \pi u(P) - \int_{\partial \Omega} u(P) \frac{\partial}{\partial n_Q} \log |P - Q| d S_Q = - f(P), \quad P \in \partial \Omega,
\end{align*}
is solved by discretizing in $\theta$ with $n$ points and applying the trapezoidal rule with $n = 2m$ (see \cite{atkinson}).  For the Bernoulli Dirichlet Ensemble (BDE), $X_m$ and $Y_m$ are Bernoulli variables taking values $\pm 1/(2m)$ with equal probability.  For the Uniform Dirichlet Ensemble (UDE), $X_m$ and $Y_m$ are uniform variables on $[-1/(2m),1/(2m)]$.

\section*{Acknowledgments}
We acknowledge the partial support of the National Science Foundation through grants NSF-DMS-130318 (TT), NSF-DMS-1300965 (PD) and NSF-DMS-07-48482 (GM) and the Australian Research Council through the Discovery Early Career Research Award (SO).  Any opinions, findings, and conclusions or recommendations expressed in this material are those of the authors and do not necessarily reflect the views of the funding sources.

\bibliographystyle{pnas}
\bibliography{/Users/trogdon/Dropbox/References/library}
%\begin{thebibliography}{10}

%\end{thebibliography}

\end{document}